\documentclass[11pt]{amsart}
\usepackage{mathrsfs}
\usepackage{amssymb}

\usepackage{titletoc}
\usepackage{amscd}
\usepackage{amsmath, amssymb}
\usepackage{amsfonts}
\usepackage[colorlinks,linkcolor=blue,citecolor=blue, pdfstartview=FitH]{hyperref}

\numberwithin{equation}{section}

\theoremstyle{plain}
\newtheorem{thm}{Theorem}[section]
\newtheorem{theorem}[thm]{Theorem}
\newtheorem{lemma}[thm]{Lemma}
\newtheorem{corollary}[thm]{Corollary}
\newtheorem{proposition}[thm]{Proposition}

\newtheorem{conjecture}[thm]{Conjecture}
\theoremstyle{definition}

\newtheorem{question}[thm]{Question}

\newtheorem{remark}[thm]{Remark}

\newtheorem{definition}[thm]{Definition}

\newtheorem{example}[thm]{Example}

\newtheorem{defn-thm}[thm]{Definition-Theorem}

\newcommand{\sE}{{\mathcal E}}

\newcommand{\sO}{{\mathcal O}}

\newcommand{\C}{{\mathbb C}}

\renewcommand{\P}{{\mathbb P}}
\newcommand{\Q}{{\mathbb Q}}
\newcommand{\R}{{\mathbb R}}

\newcommand{\qtq}[1]{\quad\mbox{#1}\quad}
\newcommand{\bp}{\bar{\partial}}
\newcommand{\Om}{\Omega}

\newcommand{\ts}{\otimes}

\newcommand{\btheorem}{\begin{theorem}}
\newcommand{\etheorem}{\end{theorem}}
\newcommand{\bproposition}{\begin{proposition}}
\newcommand{\eproposition}{\end{proposition}}
\newcommand{\bdefinition}{\begin{definition}}
\newcommand{\edefinition}{\end{definition}}
\newcommand{\bcorollary}{\begin{corollary}}
\newcommand{\ecorollary}{\end{corollary}}
\newcommand{\bproof}{\begin{proof}}
\newcommand{\eproof}{\end{proof}}
\newcommand{\bremark}{\begin{remark}}
\newcommand{\eremark}{\end{remark}}
\newcommand{\eexample}{\end{example}}
\newcommand{\bexample}{\begin{example}}
\newcommand{\la}{\langle}
\newcommand{\elemma}{\end{lemma}}
\newcommand{\blemma}{\begin{lemma}}
\newcommand{\ra}{\rangle}
\newcommand{\sq}{\sqrt{-1}}

\newcommand{\p}{\partial}

\renewcommand{\bar}{\overline}

\renewcommand{\phi}{\varphi}

\newcommand{\ee}{\end{eqnarray*}}
\newcommand{\be}{\begin{eqnarray*}}

\newcommand{\beq}{\begin{equation}}
\newcommand{\eeq}{\end{equation}}

\newcommand{\bd}{\begin{enumerate}}
\newcommand{\ed}{\end{enumerate}}

\renewcommand{\hat}{\widehat}
\renewcommand{\tilde}{\widetilde}

\renewcommand{\>}{\rightarrow}

\usepackage{fancyhdr}
\begin{document}
\title{RC-positivity,  rational connectedness and Yau's conjecture} \makeatletter
\let\uppercasenonmath\@gobble
\let\MakeUppercase\relax
\let\scshape\relax
\makeatother
\author{Xiaokui Yang}
\date{}
\address{Morningside Center of Mathematics, Academy of Mathematics and\\ Systems Science, Chinese Academy of Sciences, Beijing, 100190, China}
\address{HCMS, CEMS, NCNIS, HLM, UCAS, Academy of Mathematics and\\ Systems Science, Chinese Academy of Sciences, Beijing 100190,
China} \email{\href{mailto:xkyang@amss.ac.cn}{{xkyang@amss.ac.cn}}}

\thanks{This work
was partially supported   by China's Recruitment
 Program of Global Experts and  NSFC 11688101.}

\maketitle
\begin{abstract}
In
 this paper,   we introduce a concept of RC-positivity for  Hermitian holomorphic  vector
 bundles and prove that, if $E$ is an RC-positive vector bundle over a compact complex manifold $X$,
then  for any vector bundle $A$, there exists a positive integer
$c_A=c(A,E)$ such that
$$H^0(X,\mathrm{Sym}^{\ts \ell}E^*\ts A^{\ts k})=0$$
for $\ell\geq c_A(k+1)$ and $k\geq 0$. Moreover,  we obtain that, on
a compact K\"ahler manifold $X$,  if $\Lambda^p T_X$ is RC-positive
for every $1\leq p\leq \dim X$, then $X$ is projective and
rationally connected. As applications, we show that if a compact
K\"ahler manifold $(X,\omega)$ has positive holomorphic sectional
curvature, then $\Lambda^p T_X$ is RC-positive and
$H_{\bp}^{p,0}(X)=0$ for every $1\leq p\leq \dim X$, and in
particular, we establish that $X$ is a projective and rationally
connected manifold, which confirms a conjecture of Yau
(\cite[Problem~47]{Yau82}).
\end{abstract}



\section{Introduction}

 In
 this paper, we give a geometric interpretation of Mumford's conjecture on rational
  connectedness by using curvature conditions and  propose a differential geometric approach to attack this
  conjecture. As an application of this approach, we confirm a
  conjecture of Yau (\cite[Problem~47]{Yau82}) that a compact K\"ahler manifold with positive
  holomorphic sectional curvature  is a projective and rationally
  connected manifold.
 This project is motivated by a number of well-known conjectures proposed by Yau,
 Mumford, Demailly, Campana, Peternell and etc., and we refer to
 \cite{Yau74, Yau82, Cam92, KMM92, Kol96, DPS96, DPS96B, Pet06, HW12, BDPP13, CDP, CP14, HW15, BC15, Yang16, Cam16, CH17} and the references
 therein.\\

A projective manifold $X$ is called \emph{rationally connected} if
any two points of $X$ can be connected by some rational curve. It is
easy to show that on a rationally connected projective manifold, one
has
 $$ H^0(X,(T^*_X)^{\ts m})=0, \qtq{for every} m\geq 1.$$
A well-known conjecture of Mumford says that the converse is also
true.
\begin{conjecture}[Mumford]\label{mainconjecture}Let $X$ be a projective manifold. If $$  H^0(X,(T^*_X)^{\ts m})=0, \qtq{for
every} m\geq 1,$$ then $X$ is rationally connected.
\end{conjecture}

\noindent This conjecture holds when $\dim X\leq 3$ (\cite{KMM92})
and not much has been known in higher dimensions, and we refer to
\cite{LP17}, \cite{CDP} and \cite{Kol96} for more historical
discussions. In \cite{BC15}, Brunebarbe and Campana also proposed a
stronger conjecture that $X$ is rationally connected if and only if
\beq H^0(X, \mathrm{Sym}^{\ts \ell}\Om_X^p)=0\ \ \ \text{for every }
\ell>0\ \ \text{and}\ \ 1\leq p\leq \dim X.\label{v1}\eeq

\vskip 1\baselineskip

 In order to give geometric interpretations on rational connectedness, we introduce the following
concept for Hermitian vector bundles: \bdefinition Let $(E,h)$ be a
Hermitian holomorphic vector bundle over a complex manifold $X$ and
$R^{(E,h)}\in \Gamma(X,\Lambda^{1,1}T_X^*\ts \mathrm{End}(E))$ be
its Chern curvature tensor.  $E$ is called \emph{RC-positive} (resp.
\emph{RC-negative}) if for any nonzero local section $a\in
\Gamma(X,E)$, there exists some local section $v\in\Gamma(X,T_X)$
such that  $$ R^{(E,h)}(v,\bar v, a,\bar a)>0. \ \ \ (\qtq{resp.}
<0)
$$ \edefinition

\noindent For a line bundle $(L,h)$, it is RC-positive if and only
if its Ricci curvature has at least one positive eigenvalue at each
point of $X$. This terminology has many nice properties. For
instances, quotient bundles of RC-positive bundles are also
RC-positive; subbundles of RC-negative bundles are still RC-negative
(see Theorem \ref{mono}); the holomorphic tangent bundles of Fano
manifolds can admit RC-positive metrics (see Corollary \ref{Fano}).
The first main result of our paper is

\btheorem\label{A} Let $X$ be a \emph{compact complex manifold}. If
$E$ is an RC-positive vector bundle, then
  for any vector bundle $A$,
there exists a positive integer $c_A=c(A,E)$  such that
$$H^0(X,\mathrm{Sym}^{\ts \ell}E^*\ts A^{\ts k})=0$$
for $\ell\geq c_A(k+1)$ and $k\geq 0$.  Moreover, if $X$ is a
projective manifold, then any \emph{invertible subsheaf}
$\mathcal{F}$ of $\sO(E^*)$ is not pseudo-effective.
 \etheorem

  As a straightforward application of Theorem \ref{A} and Campana-Demailly-Peternell's  criterion
for rational connectedness (\cite[Theorem~1.1]{CDP}) (see also
\cite[Corollary~1.7]{GHS03}, \cite{Pet06},
\cite[Proposition~2.1]{LP17} and \cite[Proposition~1.3]{Cam16}), we
obtain the second main result of our paper

\btheorem\label{Z0} Let $X$ be a compact K\"ahler manifold of
complex dimension $n$. Suppose that for every $1\leq p\leq n$, there
exists a smooth Hermitian metric $h_p$ on the vector bundle
$\Lambda^p T_X$ such that $(\Lambda^p T_X, h_p)$ is RC-positive,
then $X$ is projective and rationally connected. \etheorem

\noindent  Theorem \ref{A} and Theorem \ref{Z0} also hold  if we
replace the RC-positivity by a weaker  condition defined in
Definition \ref{Def}. In the following, we shall verify that several
classical curvature conditions in differential geometry can imply
the RC-positivity in Theorem \ref{Z0}.

\bcorollary \label{B} Let $X$ be a compact K\"ahler manifold. If
there exist a Hermitian metric $\omega$ on $X$ and a (possibly
different) Hermitian metric $h$ on $T_X$ such that \beq
\mathrm{tr}_{\omega}R^{(T_X,h)} \in \Gamma(X, \mathrm{End}(T_X))\eeq
is positive definite, then $X$ is projective and rationally
connected. \ecorollary

\noindent We need to point out that, when $X$ is projective,
Corollary \ref{B} can also be implied by the ``Generalized holonomy
principle" for \emph{positive curvature} and the main theorem in
\cite{CDP}, although the precise result is not
 stated there. This special case is a refinement under the RC positive curvature
 condition. In particular, by
 the celebrated Calabi-Yau theorem (\cite{Yau78}), we obtain the
classical result of Campana(\cite{Cam92}) and Koll\'ar-Miyaoka-Mori
(\cite{KMM92}) that Fano manifolds are rationally connected.\\

 Let's describe another application of Theorem \ref{Z0}. In
his ``Problem section",  S.-T. Yau proposed the following well-known
conjecture (\cite[Problem~47]{Yau82}):

  \begin{conjecture}[Yau] \label{Yau} Let $X$ be a compact K\"ahler manifold. If $X$ has  a
  K\"ahler metric with positive holomorphic sectional curvature,
  then $X$ is a projective and rationally connected manifold.
  \end{conjecture}

  \noindent As  applications of Theorem \ref{Z0}, we confirm Yau's
  Conjecture \ref{Yau}. More generally, we
  obtain

\btheorem\label{HSCF} Let $(X,\omega)$ be a compact K\"ahler
manifold  with positive holomorphic sectional curvature.  Then for
every $1\leq p\leq \dim X$, $(\Lambda^p T_X, \Lambda^p \omega)$ is
RC-positive and $H_{\bp}^{p,0}(X)=0$. In particular, $X$ is a
projective and  rationally connected manifold. \etheorem

 \bremark
   We need to point out that,  recently, Heier-Wong also confirmed Yau's conjecture   in the special case when $X$ is
   projective (\cite{HW15}).  \emph{One should see clearly the significant difference
   in the proofs.}
 Our method crucially relies on the geometric properties of RC-positivity (Theorem \ref{A}) which \emph{we prove by using techniques in non-K\"ahler geometry},
  and also a minimum principle for K\"ahler metrics with positive holomorphic sectional curvature (Lemma \ref{linear} and Theorem \ref{key1}), while their method
 builds on an average argument and certain integration by parts on projective manifolds. \eremark

\bremark For the negative holomorphic sectional curvature case, in
the recent breakthrough  paper \cite{WY} of  Wu and Yau,   they
proved that any projective K\"ahler manifold with negative
   holomorphic sectional curvature must have ample canonical line bundle. This result was obtained by Heier et.\! al.\! earlier
    under the additional assumption of the Abundance Conjecture (\cite{HLW10}). In \cite{TY}, Tosatti and the  author proved that any
compact K\"ahler manifold with nonpositive holomorphic sectional
 curvature must have nef canonical line bundle,  and with that in hand,  we were able to drop the projectivity assumption in the
  aforementioned Wu-Yau Theorem.  More recently, Diverio and Trapani \cite{DT} further generalized the result by assuming that
   the holomorphic sectional curvature is only {\em quasi-negative,} namely, nonpositive everywhere and negative somewhere
    in the manifold. In \cite{WY1}, Wu and Yau give a direct proof of the statement that any compact K\"ahler manifold with
     quasi-negative holomorphic sectional curvature must have ample canonical line
     bundle. We refer to \cite{Wo,HLW10, WWY, HLW14, WY, WY1, TY,
     DT, YZ, T} for more details.
On the other hand, it is well-known that the anti-canonical bundle
of a compact K\"ahler manifold with positive holomorphic sectional
curvature is not necessarily nef (e.g. \cite{Yang16}). \eremark

\noindent It is clear that, the RC-positivity is defined for
Hermitian metrics. We also obtain \bcorollary\label{C} Let $X$ be a
compact K\"ahler surface. If there exists a \emph{Hermitian metric}
$\omega$ with positive holomorphic sectional curvature, then $X$ is
a projective and rationally connected manifold. Moreover, the Euler
characteristic  $\chi(X)>0$. \ecorollary

\noindent As motivated by Theorem \ref{HSCF} and Corollary \ref{C},
we propose
\begin{conjecture}\label{HSCC}  Let $X$ be a compact complex manifold of complex dimension $n>2$. Suppose
$X$ has a Hermitian metric with positive holomorphic sectional
curvature, then $\Lambda^p T_X$ admits a smooth RC-positive metric
for every $1<p<n$. In particular, if in addition $X$ is  K\"ahler,
then $X$ is a projective and rationally connected manifold.
\end{conjecture}

\noindent For some related topics on positive holomorphic sectional
curvature, we refer to \cite{HW12, HW15, ACH15, Yang16, AHZ16, YZ16,
AH17} and the references therein.\\

Organizations. In Section \ref{pre}, we recall some background
materials on Hermitian manifolds. In Section \ref{partial}, we
introduce the concept of RC-positivity for vector bundles and
investigate its geometric properties. In Section \ref{van}, we
obtain vanishing theorems for RC-positive vector bundles and prove
Theorem \ref{A}. In Section \ref{manifold}, we establish the
relation between RC-positive bundles and rational connectedness, and
prove Theorem \ref{Z0} and Corollary \ref{B}. In Section
\ref{HSC00}, we study positive holomorphic sectional curvature and
rational connectedness, and prove Theorem \ref{HSCF}, and Corollary
\ref{C}. In Section \ref{open}, we study the relations between
several open conjectures and propose some further
questions.\\

\section{Background materials}\label{pre}

Let $(E,h)$ be a Hermitian holomorphic vector bundle over a complex
manifold $X$ with Chern connection $\nabla$. Let $\{z^i\}_{i=1}^n$
be the  local holomorphic coordinates
  on $X$ and  $\{e_\alpha\}_{\alpha=1}^r$ be a local frame
 of $E$. The curvature tensor $R^{(E,h)}\in \Gamma(X,\Lambda^{1,1}T^*_X\ts \mathrm{End}(E))$ has components \beq R_{i\bar j\alpha\bar\beta}= -\frac{\p^2
h_{\alpha\bar \beta}}{\p z^i\p\bar z^j}+h^{\gamma\bar
\delta}\frac{\p h_{\alpha \bar \delta}}{\p z^i}\frac{\p
h_{\gamma\bar\beta}}{\p \bar z^j}.\eeq (Here and henceforth we
sometimes adopt the Einstein convention for summation.) We have the
trace $\mathrm{tr} R^{(E,h)}\in  \Gamma(X,\Lambda^{1,1}T^*_X)$ which
has components
$$R_{i\bar j \alpha}^\alpha=h^{\alpha\bar\beta}R_{i\bar j\alpha\bar\beta}=-\frac{\p^2\log \det(h_{\alpha\bar\beta})}{\p z^i\p\bar z^j}.$$ For any
Hermitian metric $\omega_g$ on $X$, we can also define
$\mathrm{tr}_{\omega_g} R^{(E,h)}\in \Gamma(X,\mathrm{End}(E))$
which has components
$$g^{i\bar j} R_{i\bar j \alpha}^{\beta}.$$
In particular, if $(X,\omega_g)$ is a  Hermitian manifold, then the
Hermitian vector bundle $(T_X,g)$ has  Chern curvature components
\beq R_{i\bar j k\bar \ell}=-\frac{\p^2g_{k\bar \ell}}{\p z^i\p\bar
z^j}+g^{p\bar q}\frac{\p g_{k\bar q}}{\p z^i}\frac{\p g_{p\bar
\ell}}{\p\bar z^j}.\eeq

\noindent The (first) Chern-Ricci form
$$\mathrm{Ric}(\omega_g)=\mathrm{tr}_{\omega_g} R^{(T_X,g)}
\in\Gamma(X,\Lambda^{1,1}T_X^*)$$ of $(X,\omega_g)$ has components
$$R_{i\bar j}=g^{k\bar \ell}R_{i\bar jk\bar \ell}=-\frac{\p^2\log\det(g)}{\p z^i\p\bar z^j}$$
and it is well-known that the Chern-Ricci form represents the first
Chern class of the complex manifold $X$ (up to a factor $2\pi$). The
second Chern-Ricci tensor
$$\mathrm{Ric}^{(2)}(\omega_g)=\mathrm{tr}_{\omega_g} R^{(T_X,g)}
\in\Gamma(X,\mathrm{End}(T_X))$$ of $(X,\omega_g)$ has (lowered
down) components
$$R^{(2)}_{k\bar \ell}=g^{i\bar j}R_{i\bar jk\bar \ell}.$$
If $\omega_g$ is not K\"ahler (i.e. $d\omega_g\neq 0$),
$\mathrm{Ric}(\omega_g)$ and $\mathrm{Ric}^{(2)}(\omega_g)$ are not
the same. The (Chern) scalar curvature $s_g$ of $\omega_g$ is
defined as
$$s_g=g^{i\bar j} R_{i\bar j}.$$

\bdefinition  A holomorphic vector bundle $(E,h)$ is called
\emph{Griffiths positive} if $$R_{i\bar j \alpha\bar\beta} v^i\bar
v^j a^\alpha \bar a^\beta>0$$ for any nonzero vectors $v=(v^i)$ and
$a=(a^\alpha)$.

A Hermitian (or K\"ahler) manifold  $(X,\omega)$ has positive (resp.
semi-positive) holomorphic sectional curvature, if for any nonzero
vector $\xi=(\xi^1,\cdots, \xi^n)$,
$$R_{i\bar j k\bar \ell}\xi^i\bar\xi^j\xi^k\bar\xi^\ell>0\ \ \ \text{(resp. $\geq 0$)},$$
at each point of $X$. \edefinition

\vskip 1\baselineskip

\section{RC-positive vector bundles on compact complex
manifolds}\label{partial} In the section we introduce the concept of
RC-positivity for Hermitian vector bundles and investigate its
geometric properties.

\begin{definition} Let $L$ be a holomorphic line bundle over a  compact complex
manifold $X$ with $\dim_\C X=n$.  $L$ is called \emph{$q$-positive},
if there exists a smooth Hermitian metric $h$ on $L$ such that the
Chern curvature $R^{(L,h)}=-\sq\p\bp\log h$ has at least $(n-q)$
positive eigenvalues at every point on $X$. \edefinition

\noindent When $q=n-1$, the concept of $(n-1)$ positivity has very
nice geometric interpretations. We established in
\cite[Theorem~1.5]{Yang17D} the following result:

\btheorem\label{t} Let $L$ be a line bundle over a compact complex
manifold $X$ with $\dim_\C X=n$. The following statements are
equivalent:

\bd
\item[(1)] $L$ is $(n-1)$-positive;
\item[(2)] The dual line bundle $L^{*}$ is not pseudo-effective.

\ed \etheorem

 \noindent In this paper, we extend the concept of $(n-1)$-positivity
to vector bundles.

 \bdefinition\label{Def}
A Hermitian holomorphic vector bundle $(E,h)$ over a complex
manifold $X$ is called \emph{RC-positive}  (resp.
\emph{RC-negative}) at point $q\in X$ if for any nonzero
$a=(a^1,\cdots, a^r)\in \C^r$, there exists a vector $v=(v^1,\cdots,
v^n)\in \C^n$ such that \beq \sum R_{i\bar j \alpha\bar\beta}v^i\bar
v^j a^\alpha\bar a^\beta>0\ \ \ \ (\text{resp}. <0) \eeq at point
$q$. $(E,h)$ is  called RC-positive if it is RC-positive at every
point of $X$. $E$ is called \emph{weakly RC-positive} if there
exists a smooth Hermitian metric $h$ on the tautological line bundle
$\sO_E(1)$ over $\P(E^*)$ such that $(\sO_E(1), h)$ is $(\dim
X-1)$-positive. \edefinition

 \noindent We can also define RC-semi-positivity (resp. RC-semi-negativity) in the same way. \bremark  From the definition, it
is easy to see that, \bd \item[(1)] if $(E,h)$ is RC-positive, then
for any nonzero $u\in\Gamma(X,E)$, as a Hermitian $(1,1)$-form on
$X$, $R^{(E,h)}(u,u)\in \Gamma(X, \Lambda^{1,1}T^*_X)$ has at least
one positive eigenvalue, i.e. $R^{(E,h)}(u,u)$ is $(n-1)$-positive;
\item[(2)] if a vector bundle $(E,h)$ is Griffiths positive, then
$(E,h)$ is RC-positive;
\item[(3)] if $E$ is a line bundle, then $E$ is RC-positive if and
only if $E$ is $(n-1)$-positive; \item[(4)] if $\dim X=1$, $(E,h)$
is RC-positive if and only if $(E,h)$ is Griffiths positive. \ed

\noindent  Here, we can also define RC-positivity along $k$ linearly
independent directions, i.e. for any given nonzero local section
$a\in \Gamma(X,E)$, $R^{(E,h)}(\bullet,\bullet, a, \bar a)\in
\Gamma(X,\Lambda^{1,1}T_X^*)$ is $(n-k)$-positive as a Hermitian
$(1,1)$-form on $X$. It is also a generalization of the Griffiths
positivity. This terminology  will be systematically investigated in
the forthcoming paper \cite{Yang18}.
 \eremark

\noindent By using the monotonicity formula and Theorem \ref{t}, we
obtain the following properties, which also hold for RC-positivity
along $k$ linearly independent directions.

\btheorem\label{mono} Let $(E,h)$ be a Hermitian vector bundle over
a compact complex manifold $X$. \bd

\item[(1)]  $(E,h)$ is RC-positive if and only if $(E^*,h^*)$ is
RC-negative;

 \item[(2)]  If $(E,h)$ is RC-negative, every subbundle $S$ of $E$ is RC-negative;

\item[(3)] If $(E,h)$ is RC-positive,
 every quotient bundle $Q$ of $E$ is RC-positive;

\item[(4)]  If  $(E,h)$ is RC-positive, every line subbundle $L$ of
$E^*$ is not pseudo-effective.

\item[(5)] If $(E,h)$ is an RC-positive line bundle, then   for any
pseudo-effective line bundle $L$, $E\ts L$ is RC-positive.

\ed \etheorem

\bproof $(1)$ is obvious.  $(2)$ follows from a standard
monotonicity formula.
 Let $r$ be the rank of $E$ and $s$ the rank of $S$. Without
loss of generality, we can assume,  at a fixed point $p\in X$, there
exists a local holomorphic frame $\{e_1,\cdots,e_r\}$ of $E$
centered at point $p$ such that $\{e_1,\cdots, e_s\}$ is a local
holomorphic frame of $S$. Moreover, we can assume that
$h(e_\alpha,e_\beta)(p)=\delta_{\alpha\beta}, \text{for } 1\leq
\alpha, \beta \leq r.$ Hence, the curvature tensor of $S$ at point
$p$ is \beq R^S_{i\bar j \alpha\bar\beta}=-\frac{\p^2
h_{\alpha\bar\beta}}{\p z^i\p\bar z^j}+\sum_{\gamma=1}^{s}\frac{\p
h_{\alpha\bar\gamma}}{\p z^i}\frac{\p h_{\gamma\bar\beta}}{\p\bar
z^j} \eeq where $1\leq \alpha,\beta\leq s$. The curvature tensor of
$E$ at point $p$ is \beq R^E_{i\bar j \alpha\bar\beta}=-\frac{\p^2
h_{\alpha\bar\beta}}{\p z^i\p\bar z^j}+\sum_{\gamma=1}^{r}\frac{\p
h_{\alpha\bar\gamma}}{\p z^i}\frac{\p h_{\gamma\bar\beta}}{\p\bar
z^j} \eeq where $1\leq \alpha,\beta\leq r$. It is easy to see that
\beq R^E|_{S}-R^S=\sq
\sum_{i,j}\sum_{\alpha,\beta=1}^s\left(\sum_{\gamma=s+1}^r\frac{\p
h_{\alpha\bar\gamma}}{\p z^i}\frac{\p h_{\gamma\bar\beta}}{\p\bar
z^j}\right)dz^i\wedge d\bar z^j\ts e^\alpha\ts
e^\beta.\label{cu}\eeq  Since $R^E$ is RC-negative,  for any nonzero
vector $ a_E=(a^1,\cdots, a^s, 0,\cdots, 0)\in \C^r$, there exists a
vector $v=(v^1,\cdots, v^n)$ such that
$$R^E(v,\bar v ,a_E, \bar a_E)<0.$$
Let $a=(a^1,\cdots, a^s)$. Then \beq R^S(v,\bar v, a ,\bar
a)=R^{E}(v,\bar v ,a_E, \bar a_E)-
\sum_{i,j}\sum_{\alpha,\beta=1}^s\left(\sum_{\gamma=s+1}^r\frac{\p
h_{\alpha\bar\gamma}}{\p z^i}\frac{\p h_{\gamma\bar\beta}}{\p\bar
z^j}\right)v^i\bar v^j a^\alpha\bar a^{\beta}<0.\label{keyf}\eeq The
proof of part $(3)$ is similar.\\

$(4)$. Let $L$ be a line subbundle of $E^*$ where $(E,h)$ is
RC-positive. Then by $(1)$ and $(2)$, we know $L$ is RC-negative.
Hence $L^*$ is RC-positive and $L^*$ is $(n-1)$-positive. By Theorem
\ref{t}, $L=(L^*)^*$ is not pseudo-effective.\\

$(5)$. Suppose $E\ts L$ is not RC-positive, then by Theorem \ref{t},
$E^*\ts L^*$ is pseudo-effective and so is $E^*=(E^*\ts L^*)\ts L$.
This is a contradiction since $E$ is RC-positive.
  \eproof

Let's recall some basic linear algebra. Let $V$ be a complex vector
space and $\dim_\C V=r$. Let $A\in \mathrm{End}(V)$. For any $1\leq
p\leq r$, we define $\Lambda^pA\in\mathrm{End}(\Lambda^pV)$ as \beq
\left(\Lambda^pA\right)(v_1\wedge \cdots \wedge v_p)=\sum_{i=1}^p
v_1\wedge \cdots\wedge A v_i\wedge \cdots\wedge v_p.\eeq Similarly,
for $k\geq 1$ we define $A^{\ts k}\in \mathrm{End}(V^{\ts k})$ by
\beq \left(A^{\ts k}\right)(v_1\ts\cdots \ts v_p)=\sum_{i=1}^p
v_1\ts \cdots\ts A v_i\ts \cdots\ts v_p.\eeq By choosing a basis, it
is easy to see that if $A$ is positive definite, then both
$\Lambda^pA$ and $A^{\ts k}$ are positive definite.
 We have the following important observation:

 \btheorem\label{E} Let
$(X,\omega)$ be a Hermitian manifold and $(E,h)$ be a Hermitian
holomorphic vector bundle. If \beq \mathrm{tr}_\omega R^{(E,h)}\in
\Gamma(X,\mathrm{End}(E)) \label{P}\eeq is positive definite, then
\bd

\item[(1)] $(\Lambda^pE,\Lambda^p h)$ is
RC-positive for every $1\leq p\leq \mathrm{rank}(E)$;
\item[(2)] $(E^{\ts k}, h^{\ts k})$ is RC-positive for every $k\geq
1$. \ed \etheorem

\bproof From the expression of the induced curvature tensor of
$(\Lambda^pE,\Lambda^p h)$, it is easy to see that
$$R^{(\Lambda^pE,\Lambda^p h)}=\Lambda^pR^{(E,h)}\in \Gamma(X,\Lambda^{1,1}T_X^*\ts \mathrm{End}(\Lambda^pE)),$$
and
$$\mathrm{tr}_\omega R^{(\Lambda^pE,\Lambda^p h)}=\mathrm{tr}_\omega\left(\Lambda^pR^{(E,h)}\right)=\Lambda^p\left(\mathrm{tr}_\omega R^{(E,h)}\right)\in \Gamma(X,\mathrm{End}(\Lambda^pE)),$$
is positive definite. Suppose $(\Lambda^pE,\Lambda^p h)$ is not
RC-positive, then there exists a point $q\in X$ and a nonzero local
section $a\in \Gamma(X,\Lambda^pE)$, such that for any local section
$v\in \Gamma(X,T_X)$,
$$R^{(\Lambda^pE,\Lambda^p h)}(v,\bar v, a, \bar a)\leq 0$$ at point
$q\in X$. In particular, we have
$$\left(\mathrm{tr}_\omega R^{(\Lambda^pE,\Lambda^p h)}\right)(a, \bar a)\leq
0$$ at point $q$. This is a contradiction. The proof for part $(2)$
is similar.
 \eproof

\bcorollary Let $X$ be a compact complex manifold.  If there exist a
\emph{Hermitian metric} $\omega$ on $X$ and a (possibly different)
Hermitian metric $h$ on $T_X$ such that \beq \mathrm{tr}_\omega
R^{(T_X,h)}>0\eeq then \bd
\item[(1)] $(T^{\ts k}_X,h^{\ts k})$ is  RC-positive for
every $k\geq 1$;
\item[(2)] $(\Lambda^pT^{}_X,\Lambda^p h)$ is  RC-positive for every  $1\leq p\leq \dim X$.
\ed \ecorollary

\bcorollary\label{Fano} Let $X$ be a compact K\"ahler manifold.  If
there exists a \emph{K\"ahler metric} $\omega$ such that it has
positive Ricci curvature, i.e. $$ \mathrm{Ric}(\omega)=-\sq\p\bp\log
(\omega^n)>0,
$$ then \bd \item[(1)] $(T^{\ts k}_X,\omega^{\ts k})$ is  RC-positive for every $k\geq
1$;
\item[(2)] $(\Lambda^pT^{}_X,\Lambda^p\omega^{})$ is  RC-positive for every  $1\leq p\leq \dim X$.
\ed \ecorollary

\vskip 1\baselineskip

\section{Vanishing theorems for RC-positive vector
bundles}\label{van} In this section, we derive vanishing theorems
for RC-positive vector bundles over compact complex manifolds and
prove Theorem \ref{A}.

\bproposition\label{G} Let $\sO_E(1)\>\P(E^*)$ be the tautological
line bundle of $E\>X$. Suppose $(E,h^E)$ is RC-positive, then $E$ is
weakly RC-positive over $\P(E^*)$. \eproposition \bproof The proof
follows from a standard curvature formula for $\sO_E(1)$ induced
from  $(E,h)$. Suppose $\dim_\C X=n$. Let $\pi$ be the projection
$\P(E^*)\> X$ and $L=\sO_E(1)$. Let $(e_1,\cdots, e_r)$ be the local
holomorphic frame with respect to a given trivialization on $E$ and
the dual frame on $E^*$ is denoted by $(e^1,\cdots, e^r)$. The
corresponding holomorphic coordinates on $E^*$ are denoted by
$(W_1,\cdots, W_r)$. There is a local section $e_{L^*}$ of $L^*$
defined by \beq e_{L^*}=\sum_{\alpha=1}^r W_\alpha e^\alpha\eeq Its
dual section is denoted by $e_L$. Let $h^L$ the induced quotient
metric on $L$ by the morphism $(\pi^*E,\pi^*h^E)\>L$. If
$\left(h_{\alpha\bar\beta}\right)$ is the matrix representation of
$h$ with respect to the basis $\{e_\alpha\}_{\alpha=1}^r$, then
$h^L$ can be written as \beq
h^L=\frac{1}{h^{L^*}(e_{L^*},e_{L^*})}=\frac{1}{\sum
h^{\alpha\bar\beta}W_\alpha\bar W_\beta}\label{inducedmetric} \eeq
 The curvature of $(L,h^L)$ is \beq
R^{h^L}=-\sq\p\bp\log h^L=\sq\p\bp\log\left(\sum
h^{\alpha\bar\beta}W_\alpha\bar W_\beta\right)
\label{inducedcurvature}\eeq where $\p$ and $\bp$ are operators on
the total space $\P(E^*)$. We fix a point $p\in \P(E^*)$, then there
exist local holomorphic coordinates
 $(z^1,\cdots, z^n)$ centered at point $s=\pi(p)$ and local holomorphic basis $\{e_1,\cdots, e_r\}$ of $E$ around $s$ such that
 \beq h_{\alpha\bar\beta}=\delta_{\alpha\bar\beta}-R_{i\bar j \alpha\bar\beta}z^i\bar z^j+O(|z|^3) \label{normal}\eeq
Without loss of generality, we assume $p$  is the point $(0,\cdots,
0,[a_1,\cdots, a_r])$ with $a_r=1$. On the chart $U=\{W_r=1\}$ of
the fiber $\P^{r-1}$, we set $w^A=W_A$ for $A=1,\cdots, r-1$. By
formula (\ref{inducedcurvature}) and (\ref{normal}) \beq
R^{h^L}(p)=\sq\left(\sum R_{i\bar j\alpha\bar\beta}\frac{a_\beta
\bar a_\alpha}{|a|^2}dz^i\wedge d\bar
z^j+\sum_{A,B=1}^{r-1}\left(1-\frac{a_B\bar
a_A}{|a|^2}\right)dw^A\wedge d\bar w^B\right) \eeq where
$|a|^2=\sum\limits_{\alpha=1}^r|a_\alpha|^2$. If $R^{E}$  is
RC-positive, then the $(n\times n)$ Hermitian matrix \beq \left[\sum
R_{i\bar j\alpha\bar\beta}\frac{a_\beta \bar
a_\alpha}{|a|^2}\right]_{i,j=1}^n\eeq has at least one positive
eigenvalues. In particular, the curvature $R^{h^L}$ of $(L,h^L)$ has
at least $r$ positive eigenvalues over the projective bundle
$\P(E^*)$. Since $\dim \P(E^*)=n+r-1$, we know $L$ is
$(n-1)$-positive. \eproof

\noindent We also have the following observation for RC-positive
vector bundles.
 \bproposition\label{kahlerprojective}
Let $X$ be a compact complex manifold. If $E$ is RC-positive, then
\beq H^0(X,E^*)=0.\eeq

\eproposition

\bproof Let  $(E, h)$ be RC-positive and $\sigma\in H^0(X,E^*)$. We
have \beq \p\bp|s|^2_g=\la\nabla s, \nabla s\ra_g-R^{E^*}(\bullet,
\bullet, s,\bar s).\label{max}\eeq Suppose $|s|^2_g$ attains its
maximum at some point $p$ and $|s|^2_g(p)>0$. Since $(E,h)$ is
RC-positive, the induced bundle $(E^*,g)$ is RC negative, i.e., at
point $q$, there exists a nonzero vector $v$ such that
$$R^{E^*}(v,
\bar v, s,\bar s)<0.$$ By applying maximum principle to (\ref{max}),
we get a contradiction. Hence, we deduce $s=0$ and $H^0(X,E^*)=0$.
\eproof

\bcorollary\label{kp} Let $X$ be a compact K\"ahler manifold.
Suppose $\Lambda^2T_X$ is $RC$-positive, then $X$ is projective.
\ecorollary

\bproof By Proposition \ref{kahlerprojective}, we have
$H^{2,0}_{\bp}(X)=H^{0,2}_{\bp}(X)=0$. Hence, by the Kodaira theorem
(\cite[Theorem~1]{kod54}, see also \cite[Proposition~ 3.3.2 and
Corollary~5.3.3]{Huy05}), the K\"ahler manifold $X$ is projective.
\eproof

\blemma \label{A1} Let $X$ be a compact complex manifold of complex
dimension $n$. If $E$ is weakly RC-positive, then for any vector
bundle $A$ over $\P(E^*)$, there exists a positive integer
$c_A=c(A,E)$ such that \beq H^{p,q}(\P(E^*), \sO_E(\ell)\ts A^{\ts
k})=0 \label{vanishing} \eeq when $\ell\geq c_A(k+1) $, $k\geq 0$,
$p\geq 0$ and $q>n-1$.  \elemma

\bproof
 The proof follows from  an
Andreotti-Grauert type  vanishing theorem. Let $F=\Om^p_{\P(E^*)}\ts
K^{-1}_{\P(E^*)}$, then we have \beq H^{p,q}(\P(E^*), \sO_E(\ell)\ts
A^{\ts k})\cong H^{n+r-1,q}(\P(E^*), \sO_E(\ell)\ts A^{\ts k}\ts
F).\eeq  Suppose $E$ is weakly RC-positive, i.e., there exists a
smooth Hermitian metric on $\sO_E(1)$ such that its curvature has at
least $r$-positive eigenvalues at each point where
$r=\mathrm{rank}(E)$.  As shown in \cite[Proposition~2.2]{Yang17D},
by a conformal change of the background metric,  there exists  a
smooth Hermitian metric $\omega$ on $\P(E^*)$ such that $(\sO_E(1),
h)$ is uniformly $(n-1)$-positive. That is, in local holomorphic
coordinates of $\P(E^*)$, at some point $p\in X$,
$$\omega=\sq \sum_i dz^i\wedge d\bar z^i, \ \ \ \ R^{(\sO_E(1),h)}=\sq\sum_i \lambda_i dz^i\wedge d\bar z^i$$
where $\lambda_1\geq \cdots \geq \lambda_{n+r-1}$,  $\lambda_r>0$
and  \beq \lambda_r+\cdots+\lambda_{n+r-1}>0.\eeq Let
$\sE=\sO_E(\ell)\ts A^{\ts k}\ts F$ and $h^{\sE}$ be the induced
metric on $\sE$ where we fix an arbitrary smooth metric $h^A$ on
$A$. Without loss of generality, we assume $A$ is a line bundle. By
standard Bochner formulas on compact complex manifolds
(\cite[Theorem~1.4 in Chapter VII]{Dem}), one has  \beq
\Delta_{\bp_\sE}=\hat{\Delta}_{\p_\sE}+\left[R^{\sE},\Lambda_\omega\right]+T_\omega\eeq
where
 $\tau=[\Lambda_\omega,
\p\omega]$, $\nabla^\sE=\p_\sE+\bp_{\sE}$ is the decomposition of
the Chern connection, $
T_\omega=[\Lambda_\omega,[\Lambda_\omega,\sq\p\bp\omega]]-[\p\omega,(\p\omega)^*]$
and
$$\hat{\Delta}_{\p_\sE}=(\p_\sE+\tau)(\p_\sE^*+\tau^*)+(\p_\sE^*+\tau^*)(\p_\sE+\tau).$$
Hence, for any $s\in H^{n+r-1,q}(\P(E^*),\sE)$, we have
$$\left\la [R^{\sE},\Lambda_\omega]s,s \right\ra =\left\la\ell \left[ R^{(\sO_E(1),h)},\Lambda_\omega\right]s+k \left[R^{(A,h^A)},\Lambda_\omega\right]s+ \left[R^{(F,h^F)},\Lambda_\omega\right]s,s \right\ra.$$ It
is obvious that there exist uniform constants $c_0\in \R$ and
$c_1\in \R$ such that $$ \left\la
\left[R^{(A,h^A)},\Lambda_\omega\right]s,s \right\ra\geq c_0|s|^2
\qtq{and }\la T_\omega s, s\ra +\left\la
\left[R^{(F,h^F)},\Lambda_\omega\right]s,s \right\ra\geq c_1|s|^2$$
hold on $\P(E^*)$. On the other hand, a straightforward computation
shows that when $q>n-1$
$$\ell\left\la \left[ R^{(\sO_E(1),h)},\Lambda_\omega\right]s,s \right\ra\geq \ell\left(\lambda_r+\cdots+\lambda_{n+r-1}\right)|s|^2. $$
We set $$\gamma=\inf_{p\in \P(E^*)}
\left(\lambda_r+\cdots+\lambda_{n+r-1}\right)>0$$ and pick a
positive number $c_A=c(A,E)$ such that
$$c_A\cdot \gamma +c_0>0,\ \ c_A\cdot \gamma+c_1>1.$$
Hence, for any $s\in H^{n+r-1,q}(\P(E^*),\sE)$, $q>n-1$ and
$\ell\geq c_A(k+1)$, we have \be
0=\left\la\Delta_{\bp_\sE}s,s\right\ra&=&\left\la\hat{\Delta}_{\p_\sE}s,s\right\ra+\left\la\left[R^{\sE},\Lambda_\omega\right]s,s\right\ra+\left\la
T_\omega s,s\right\ra\\
&\geq & \left\la\hat{\Delta}_{\p_\sE}s,s\right\ra +k(c_A\cdot\gamma
+c_0)|s|^2+(c_A\cdot \gamma+c_1)|s|^2\\
&\geq & \left\la\hat{\Delta}_{\p_\sE}s,s\right\ra+|s|^2. \ee An
integration over $\P(E^*)$ shows $s=0$.
 \eproof

\btheorem\label{A2}  Let $E$ be a weakly RC-positive vector bundle
over a compact complex manifold $X$. Then
  for every vector bundle $A$,
there exists a positive integer $c_A=c(A,E)$  such that \beq
H^0(X,\mathrm{Sym}^{\ts \ell}E^*\ts A^{\ts k})=0\label{00} \eeq when
$\ell\geq c_A(k+1)$ and $k\geq 0$. In particular, (\ref{00}) holds
if $E$ is RC-positive.
 \etheorem
\bproof By the classical Le Potier isomorphism over compact complex
manifolds (e.g. \cite[Theorem~5.16]{SS}), we have
$$H^{q}\left(\P(E^*), \sO_E(\ell)\ts \Om_{\P(E^*)}^p\ts (\pi^*A^*)^{\ts k}\right)=H^{q}\left(X,\mathrm{Sym}^{\ts \ell }E \ts \Om_X^p \ts (A^*)^{\ts k}\right)$$
where $\pi:\P(E^*)\>X$ is the projection. By Lemma \ref{A1}, if we
take $p=q=n$,
$$H^{n,n}\left(\P(E^*), \sO_E(\ell)\ts (\pi^*A^*\right)^{\ts k})=H^{n,n}\left(X,\mathrm{Sym}^{\ts \ell }E \ts (A^*)^{\ts k}\right)=0$$
when $\ell\geq c_A(k+1)$ and $k\geq 0$. By the Serre duality, we
obtain (\ref{00}). \eproof

\vskip 1\baselineskip

\noindent \emph{The proof of Theorem \ref{A}.}   The first part of
Theorem \ref{A} is contained in  Theorem \ref{A2}. For the second
part,
 suppose-to the contrary-that there
exists a pseudo-effective invertible sheaf $\mathcal{F}\subset
\sO(E^*)$. It is well-known that there exists a very ample line
bundle $A$ on the projective manifold such that
$$H^0(X,\mathcal{F}^{\ts \ell}\ts A)\neq 0$$
for all $\ell\geq 0$, and so
$$H^0(X,\mathrm{Sym}^{\ts \ell}E^*\ts A)\neq 0.$$
This is a contradiction. \qed

\vskip 2\baselineskip

\section{RC positivity and rational connectedness}\label{manifold}
In this section, we prove Theorem \ref{Z0} and Corollary \ref{B}.
The proofs rely on Theorem \ref{A},  Corollary \ref{kp} and a
classical criterion for
 rational connectedness proved in \cite[Theorem~1.1]{CDP} ( see also \cite{Pet06}, \cite[Corollary~1.7]{GHS03}, \cite[Proposition~2.1]{LP17} and \cite[Proposition~1.3]{Cam16}.)\\

\noindent \emph{The proof of Theorem \ref{Z0}}.  By Corollary
\ref{kp}, we know $X$ is a projective manifold. The rational
connectedness of $X$  follows from Theorem \ref{A} and
\cite[Theorem~1.1]{CDP} (or \cite[Proposition~1.3]{Cam16},
\cite{Pet06}). For readers' convenience, we give a proof here
following \cite[Theorem~1.1]{CDP}. Indeed, if $(\Lambda^p T_X, h_p)$
is RC-positive for every $1\leq p\leq \dim X$, then by Theorem
\ref{A}, any invertible subsheaf $\mathcal{F}$ of $\Om_X^p$ can not
be pseudo-effective.   In particular, when $p=n$, $K_X$ is not
pseudo-effective. Thanks to \cite{BDPP13}, $X$ is uniruled. Let
$\pi:X\dashrightarrow Z$ be the associated MRC fibration of $X$.
After possibly resolving the singularities of $\pi$ and $Z$, we may
assume that $\pi$ is a proper morphism and $Z$ is smooth. By a
result of Graber, Harris and Starr \cite[Corollary~1.4]{GHS03}, it
follows that the target $Z$ of the MRC fibration is either a point
or a positive dimensional variety which is not unruled. Suppose $X$
is not rationally connected, then $\dim Z\geq 1$. Hence $Z$ is not
uniruled, by \cite{BDPP13} again, $K_Z$ is pseudo-effective. Since
 $K_Z=\Om^{\dim Z}_Z\subset
\Om_X^{\dim Z}$ is pseudo-effective, we get a contradiction. Hence
$X$ is rationally connected. \qed

\vskip 1\baselineskip

\noindent \emph{The proof of Corollary \ref{B}.} It follows Theorem
\ref{E} and Theorem \ref{Z0}. If the trace $\mathrm{tr}_{\omega}
R^{(T_X,h)}$ is positive definite, by Theorem \ref{E}, $(\Lambda^p
T_X, \Lambda^p h)$ is RC-positive for every $1\leq p\leq \dim X$.
Hence, by Theorem \ref{Z0}, $X$ is projective and rationally
connected. \qed

\vskip 2\baselineskip

\section{A proof of Yau's conjecture on positive holomorphic sectional curvature}\label{HSC00}

In this section, we describe more applications of Theorem \ref{Z0}
on RC-positive vector bundles and  prove Theorem \ref{HSCF} and
Corollary \ref{C}. More precisely, we confirm Yau's conjecture that
if a compact K\"ahler manifold has positive holomorphic sectional
curvature, then it is projective and rationally connected.\\

\noindent \textbf{6.1. Compact K\"ahler manifolds with positive
holomorphic sectional curvature.} We begin with an algebraic
curvature relation on a compact K\"ahler manifold $(X,\omega)$. By
the K\"ahler symmetry, we have
$$R_{i\bar j k\bar\ell}=R_{k\bar\ell i \bar j}=R_{k\bar j i\bar \ell}.$$
 At a given point $q\in X$, the minimum holomorphic
sectional curvature is defined as
$$\min_{W\in T_qX,|W|=1}H(W),$$ where $H(W):=R(W,\bar W,
W,\bar W)$. Since $X$ is of finite dimension, the minimum can be
attained. The following result is essentially obtained in
\cite[Lemma~4.1]{Yang17ann} (see also some variants in
\cite[p.~312]{Goldberg}, \cite[p.~136]{Brendle} and
\cite[Lemma~1.4]{BKT}). For the sake of completeness, we include a
proof here.

\blemma\label{linear} Let $(X,\omega)$ be a compact K\"ahler
manifold and $q$ be an arbitrary point on $X$. Let $e_1\in T_qX$ be
a unit vector which minimizes the holomorphic sectional curvature
$H$ of $\omega$ at point $q$, then \beq 2R(e_1,\bar e_1,W,\bar
W)\geq \left(1+|\la W,e_1\ra|^2\right) R(e_1,\bar e_1,e_1,\bar e_1)
\eeq for every unit vector $W\in T_qX$.

 \elemma \bproof
 Let $e_2\in T_qX$ be any unit vector  orthogonal
to $e_1$. Let $$f_1(\theta)= H(\cos(\theta )e_1+\sin(\theta) e_2),\
\ \ \ \theta\in \R.$$ Then we have the expansion \be
f_1(\theta)&=&\cos^4(\theta)R_{1\bar 1 1\bar
1}+\sin^4(\theta)R_{2\bar 2 2 \bar
2}\\&&+2\sin(\theta)\cos^3(\theta)\left[R_{1\bar 1 1\bar 2}+R_{2\bar 1 1\bar 1}\right]+2\cos(\theta)\sin^3(\theta)\left[R_{1\bar 2 2\bar 2}+R_{2\bar 1 2\bar 2}\right]\\
&&+\sin^2(\theta)\cos^2(\theta)[4R_{1\bar 1 2\bar 2}+R_{1\bar 21\bar
2}+R_{2\bar 1 2\bar 1}].\ee Since $f_1(\theta)\geq R_{1\bar 1 1\bar
1}$ for all $\theta\in \R$ and  $f_1(0)=R_{1\bar 11\bar 1}$, we have
$$f_1'(0)=0\qtq{ and} f_1''(0)\geq 0.$$ By a straightforward
computation, we obtain \beq f_1'(0)=2(R_{1\bar 1 1\bar 2}+R_{2\bar 1
1\bar 1})=0,\ \ f_1''(0)=2\left(4R_{1\bar 1 2\bar 2}+R_{1\bar 21\bar
2}+R_{2\bar 1 2\bar 1}\right)- 4R_{1\bar 1 1\bar 1}\geq 0.
\label{a}\eeq Similarly, if we set $f_2(\theta)=H(\cos(\theta
)e_1+\sq \sin(\theta) e_2)$, then \be
f_2(\theta)&=&\cos^4(\theta)R_{1\bar 1 1\bar
1}+\sin^4(\theta)R_{2\bar 2 2 \bar
2}\\&&+2\sq \sin(\theta)\cos^3(\theta)\left[- R_{1\bar 1 1\bar 2}+ R_{2\bar 1 1\bar 1}\right]+2\sq\cos(\theta)\sin^3(\theta)\left[- R_{1\bar 2 2\bar 2}+ R_{2\bar 1 2\bar 2}\right]\\
&&+\sin^2(\theta)\cos^2(\theta)[4R_{1\bar 1 2\bar 2}-R_{1\bar 21\bar
2}-R_{2\bar 1 2\bar 1}].\ee

\noindent  From  $f'_2(0)=0$ and $f''_2(0)\geq 0$, one can see \beq
-R_{1\bar 1 1\bar 2}+R_{2\bar 1 1\bar 1}=0,\ \ \ \ \ \ \
2\left(4R_{1\bar 1 2\bar 2}-R_{1\bar 21\bar 2}-R_{2\bar 1 2\bar
1}\right)-4R_{1\bar 1 1\bar 1}\geq 0.\label{b}\eeq Hence,  from
(\ref{a}) and (\ref{b}), we obtain \beq R_{1\bar 1 1\bar 2}=R_{1\bar
1 2\bar 1}=0, \qtq{and} 2R_{1\bar 1 2\bar 2}\geq R_{1\bar 1 1 \bar
1}.\label{bbb}\eeq For an arbitrary unit vector $W\in T_qX$, if $W$
is parallel to $e_1$, i.e. $W=\lambda e_1$ with $|\lambda|=1$,
$$2R(e_1,\bar e_1, W,\bar W)=2R(e_1,\bar e_1, e_1, \bar e_1).$$
Suppose $W$ is not parallel to $e_1$. Let $e_2$ be the unit vector
$$e_2=\frac{W-\la W,e_1\ra e_1}{|W-\la W,e_1\ra e_1|}.$$ Then $e_2$
is a unit vector orthogonal to $e_1$ and
$$W=ae_1+be_2$$ where  $a=\la W,e_1\ra$,  $b=|W-\la W,e_1\ra e_1|$ and
$|a|^2+|b|^2=1$. Hence $$2R(e_1,\bar{e}_1,W,\bar W)=2|a|^2R_{1\bar 1
1\bar 1}+2|b|^2 R_{1\bar 12\bar 2},$$ since we have $R_{1\bar 1
1\bar 2}=R_{1\bar 1 2\bar 1}=0$ by formula (\ref{bbb}). By formula
(\ref{bbb}) again,
$$2R(e_1,\bar{e}_1,W,\bar W)\geq (2|a|^2+|b|^2)R_{1\bar 1 1\bar 1}=(1+|a|^2) R_{1\bar 1 1\bar 1}$$
which completes the proof of Lemma \ref{linear}.
 \eproof

\bremark The K\"ahler condition is substantially used in the proof
of Lemma \ref{linear}. \eremark

\noindent Now we are ready to prove Theorem \ref{HSCF}, that is
\btheorem\label{key1} Let $(X,\omega)$ be a compact K\"ahler
manifold with positive holomorphic sectional curvature.  Then for
every $1\leq p\leq \dim X$, $(\Lambda^p T_X, \Lambda^p \omega)$ is
RC-positive and $H_{\bp}^{p,0}(X)=0$. In particular, $X$ is a
projective and rationally connected manifold. \etheorem

\bproof Suppose  $(\Lambda^p T_X, \Lambda^p \omega)$ is not
RC-positive. By  Definition \ref{Def}, there exist a point $q\in X$
and a nonzero vector $a\in \Gamma(X,E)$ where $E=\Lambda^p T_X$ such
that the Hermitian $(1,1)$-form \beq R^{E}(\bullet, \bullet, a, \bar
a)\in \Gamma(X,\Lambda^{1,1}T_X^*)\label{ke}\eeq is semi-negative at
point $q$. We choose $e_1\in \Gamma(X,T_X)$ at point $q$ such that
$$R(e_1,\bar e_1,e_1, \bar e_1)=H(e_1)=\min_{W\in T_qX,|W|=1}H(W)>0.$$
Hence, by Lemma \ref{linear}, \beq 2R(e_1,\bar e_1,W,\bar W)\geq
\left(1+|\la W,e_1\ra|^2\right) R(e_1,\bar e_1,e_1,\bar e_1)>0 \eeq
for every unit vector $W\in T_qX$. In particular,
$$R(e_1,\bar e_1, \bullet, \bullet)\in \Gamma(X, \mathrm{End}(T_X))$$
is positive definite at point $q$. Hence,\beq  R^{E}(e_1,\bar e_1,
\bullet, \bullet)\in \Gamma(X, \mathrm{End}(\Lambda^pT_X))\eeq is
also positive definite at point $q$. Therefore,
$$R^{E}(e_1,\bar e_1,a,\bar a)>0$$
which is a contradiction to (\ref{ke}). Hence, we deduce that
$(\Lambda^p T_X, \Lambda^p \omega)$ is RC-positive. By Theorem
\ref{Z0}, $X$ is projective and rationally connected.
 \eproof

\noindent\textbf{6.2. Compact complex manifolds with non-negative
holomorphic sectional curvature.} In this subsection, we investigate
Hermitian metrics with non-negative holomorphic sectional curvature.
\bproposition\label{S} Let $(X,\omega)$ be a compact Hermitian
manifold with semi-positive holomorphic sectional curvature. If the
holomorphic sectional curvature is not identically zero, then

\bd \item[(1)]  there exists a Gauduchon metric $\omega_G$ on $X$
such that
$$\int_X \mathrm{Ric}(\omega_G)\wedge \omega_G^{n-1}>0;$$
\item[(2)] $K_X$ is not pseudo-effective;

\item[(3)] there exists a  Hermitian metric $h$ on $K_X^{-1}$ such that $(K_X^{-1}, h)$ is
RC-positive.

\ed

\noindent Moreover, if in addition, $X$ is projective, then $X$ is
uniruled.

 \bremark  Proposition \ref{S} is a straightforward application of
Theorem \cite[Theorem~1.2]{Yang16}, \cite[Theorem~4.1]{Yang17D} and
the classical result of
 \cite{BDPP13}. To demonstrate the essential difficulty in proving
 Conjecture \ref{HSCC} for higher dimensional compact complex manifolds and also the significant difference from the K\"ahler case, we
 include
 a detailed proof for Proposition \ref{S}.
  \eremark

 \bproof  By \cite[Theorem~4.1]{Yang17D}, we
 know $(1)$, $(2)$ and $(3)$ are mutually equivalent. Hence, we only
 need to prove one of them, for instance  $(1)$. We follow the steps in \cite[Theorem~4.1]{Yang16} for readers' convenience. At a given point $p\in X$, the maximum holomorphic
sectional curvature is defined to be
$$ H_p:=\max_{W\in T_pX,|W|=1}H(W),$$ where $H(W):=R(W,\bar W,
W,\bar W)$.  Suppose the holomorphic sectional curvature is not
identically zero, i.e.  $ H_p>0$ for some $p\in X$. For any $q\in
X$. We assume $g_{i\bar j}(q)=\delta_{ij}$. If $\dim_{\C} X=n$ and
$[\xi^1,\cdots,\xi^n]$ are the homogeneous coordinates on
$\P^{n-1}$, and $\omega_{FS}$ is the Fubini-Study metric of
$\P^{n-1}$. At point $q$, we have the following well-known identity:
\beq\int_{\P^{n-1}}R_{i\bar j k\bar \ell}\frac{\xi^i\bar
\xi^j\xi^k\bar\xi^\ell}{|\xi|^4}\omega^{n-1}_{FS}=R_{i\bar j k\bar
\ell}\cdot
\frac{\delta_{ij}\delta_{k\ell}+\delta_{i\ell}\delta_{kj}}{n(n+1)}=\frac{s+\hat
s}{n(n+1)}.\label{scalar}\eeq where $ s$ is the Chern scalar
curvature of $\omega$ and $\hat s$ is defined as \beq \hat
s=g^{i\bar \ell} g^{k\bar j}R_{i\bar j k\bar \ell}.\eeq Hence if
$(X,\omega)$ has semi-positive holomorphic sectional curvature, then
$s+\hat s$ is a non-negative function on $X$. On the other hand, at
point $p\in X$,  $s+\hat s$ is strictly positive since $H_p>0$. By
(\ref{scalar}), the integrand is quasi-positive over $\P^{n-1}$, and
so $s+\hat s$ is strictly positive at $p\in X$. By \cite[Section
~4]{LY14},  we have the relation \beq s=\hat s+\la
\bp\bp^*\omega,\omega\ra.\label{2}\eeq  Therefore, we have \beq
\int_X \hat s \omega^n=\int_X s\omega^n-\int_X
|\bp^*\omega|^2\omega^n.\label{3}\eeq
 Let
$\omega_G=f_0^{\frac{1}{n-1}}\omega$ be a Gauduchon metric ( i.e.
$\p\bp\omega_G^{n-1}=0$ ) in the conformal class of $\omega$ for
some  positive weight function $f_0\in C^\infty(X)$. Let $s_G,\hat
s_G$ be the corresponding scalar curvatures with respect to the
Gauduchon metric $\omega_G$. Then we have
\begin{eqnarray} \int_Xs_G\omega_G^n\nonumber&=&-n\int_X \sq
\p\bp\log\det(\omega_G)\wedge \omega_G^{n-1}\\
\nonumber&=&-n\int_X f_0\sq\p\bp\log\det(\omega)\wedge
\omega^{n-1}\\
&=&\int_X f_0 s\omega^n.\label{4} \end{eqnarray}
 By using a similar equation as (\ref{3}) for $s_G,\hat s_G$ and $\omega_G$, we obtain \beq  \label{6}\int_X
\hat s_G\omega_G^n=\int_X f_0 \hat s\omega^n.\eeq
 Therefore, if $s+\hat s$ is quasi-positive,
we obtain \begin{eqnarray} \int_X
s_G\omega_G^n\nonumber&=&\frac{\int_X (s_G+\hat
s_G)\omega_G^n}{2}+\frac{\int_X (s_G-\hat
s_G)\omega_G^n}{2}\\&=&\frac{\int_X (s_G+\hat
s_G)\omega_G^n}{2}+\frac{\|\bp^*_G\omega_G\|^2}{2}=\frac{\int_X
f_0(s+\hat s)\omega^n}{2}+\frac{\|\bp^*_G\omega_G\|^2}{2}>0
\label{balanced}\end{eqnarray} where the third equation follows from
(\ref{4}) and (\ref{6}). \eproof

\eproposition

\vskip 1\baselineskip

\bcorollary\label{HSC} Let $X$ be a compact complex manifold.
Suppose $X$ has a Hermitian metric with positive holomorphic
sectional curvature, then \bd \item[(1)] $T_X$ is RC-positive;
\item[(2)] $K_X^{-1}=\det T_X$ is RC-positive.
\ed In particular, if in addition, $X$ is projective, then $X$ is
uniruled. \ecorollary

\bproof  From the definition of positive holomorphic holomorphic
sectional curvature and RC-positivity, it is easy to see that if a
Hermitian metric $\omega$ has positive holomorphic sectional
curvature, then $(T_X,\omega)$ is RC-positive. On the other hand, by
Proposition \ref{S},  if $\omega$ has positive holomorphic sectional
curvature, then there exists a (possibly different) Hermitian metric
$\tilde h$ on $K_X^{-1}$ such that $(K_X^{-1}, \tilde h)$ is
RC-positive.\eproof

\vskip 1\baselineskip

\noindent\emph{The proof of Corollary \ref{C}.} We first show $X$ is
indeed projective. Suppose $\sigma\in H^0(X,K_X)$ is not zero, then
$K_X$ is $\Q$-effective. However, by Proposition \ref{S}, $K_X$ is
not pseudo-effective. This is a contradiction. Hence we deduce
$H_{\bp}^{2,0}(X)=H_{\bp}^{0,2}(X)=H^0(X,K_X)=0$. We know $X$ is
projective.
 Now Corollary \ref{C} follows from
Theorem \ref{Z0} and Corollary \ref{HSC}. \qed

\vskip 1\baselineskip

Let $f:X\>Y$ be a smooth submersion between projective manifolds. It
is well-known that if $X$ is rationally connected, then $Y$ is
rationally connected. As  analogous to this result, we have:

\bproposition Let $f:X\>Y$ be a smooth submersion between compact
complex manifolds. Suppose $X$ admits a Hermitian metric $h$ such
that $(\Lambda^p T_X,  h)$ is RC-positive for some $1\leq p\leq \dim
Y$, then  $\Lambda^p T_Y$ admits an RC-positive Hermitian metric.
\eproposition \bproof It follows from part $(3)$ of Theorem
\ref{mono}. \eproof

\noindent Similarly, we have

\bcorollary Let $f:X\>Y$ be a smooth submersion between compact
complex manifolds. Suppose $X$ admits a Hermitian metric $h$ with
positive holomorphic sectional curvature, then $Y$ has a Hermitian
metric with positive holomorphic sectional curvature. In particular,
if in addition,  $Y$ is projective, then $Y$ is uniruled.
\ecorollary

\bproof It follows from part $(3)$ of Theorem \ref{mono}, formula
(\ref{keyf}) and Theorem \ref{HSC}. \eproof

\vskip 1\baselineskip

\section{RC positivity and Mumford's conjecture}\label{open}

In this section, we gather several conjectures in complex algebraic
geometry and  give their differential geometric interpretations.\\

\noindent \textbf{7.1. Mumford's conjecture and uniruledness
conjecture} In \cite[Theorem~4.1]{Yang17D} and
\cite[Corollary~1.6]{Yang17D}, we proved that

\btheorem\label{AA} Let $L$ be a line bundle over a compact complex
manifold $X$. The following are equivalent: \bd
\item[(1)] the dual line bundle $L^*$ is not pseudo-effective;
\item[(2)] $L$ is RC-positive;
\item[(3)] there exist a smooth Hermitian metric $h$ on $L$ and a smooth Hermitian metric $\omega$ on $X$ such that the scalar
curvature  $\mathrm{tr}_\omega (-\sq\p\bp\log h)>0$. \ed Moreover,
if $X$ is projective, then they are also equivalent to \bd\item[(4)]
  for any ample line bundle $A$,
there exists a positive integer $c_A$  such that
$$H^0(X, (L^*)^{\ts \ell}\ts A^{\ts k})=0$$
for $\ell\geq c_A(k+1)$ and $k\geq 0$.\ed\etheorem

\noindent The classical result of \cite{BDPP13} says that a
projective manifold is uniruled if and only if the canonical bundle
$K_X$ is not pseudo-effective. Hence, one can formulate the
uniruledness conjecture as

\begin{conjecture}\label{uniruledc} Let $X$ is a projective manifold. Then $\kappa(X)=-\infty$ is equivalent to one (and hence
all) of the following \bd

\item[(1)] $X$ is uniruled.
\item[(2)] $K_X$ is not pseudo-effective;
\item[(3)] $K_X^{-1}$ is RC-positive;
\item[(4)] there exists a Hermitian metric $\omega$ on $X$ with positive (Chern) scalar
curvature; \item[(5)]
  for any ample line bundle $A$,
there exists a positive integer $c_A$  such that, for $\ell\geq
c_A(k+1)$ and $k\geq 0$
$$H^0(X, K_X^{\ts \ell}\ts A^{\ts k})=0.$$\ed

\end{conjecture}

\begin{conjecture}[Mumford]\label{mainconjecture1}Let $X$ be a projective manifold. If $$  H^0(X,(T^*_X)^{\ts m})=0, \qtq{for
every} m\geq 1,$$ then $X$ is rationally connected.
\end{conjecture}

\noindent It is well-known that the uniruledness conjecture can
imply Conjecture \ref{mainconjecture1} (e.g.
\cite[Corollary~1.7]{GHS03}, see also Proposition \ref{eq}).

\begin{conjecture}\label{geo1} Let $X$ be a projective manifold. If  $$  H^0(X,(T^*_X)^{\ts m})=0, \qtq{for
every} m\geq 1,$$ then one (and hence all) of the following holds

\bd \item[(1)] $X$ is uniruled;

\item[(2)] $K_X$ is not pseudo-effective;
\item[(3)] $K_X^{-1}$ is RC-positive;
\item[(4)] there exists a
Hermitian metric $\omega$ on $X$ with positive (Chern) scalar
curvature;

\item[(5)]
  for any ample line bundle $A$,
there exists a positive integer $c_A$  such that, for $\ell\geq
c_A(k+1)$ and $k\geq 0$
$$H^0(X, K_X^{\ts \ell}\ts A^{\ts k})=0.$$

\ed

\end{conjecture}

\bproposition\label{eq} We have the following relations \beq
\qtq{Conjecture} \ref{uniruledc}\ \Longrightarrow \qtq{Conjecture}
\ref{mainconjecture1}\ \Longleftrightarrow \qtq{Conjecture}
\ref{geo1}.\eeq \eproposition

\bproof Conjecture \ref{uniruledc} $\Longrightarrow$ Conjecture
\ref{geo1}. Since $K_X=\det T_X^*$, it is well-known that for any
positive integer $\ell$, $K_X^{\ts \ell}$ is a subbundle of
$(T^*_X)^{\ts m}$ for some large $m$. Hence $H^0(X,(T^*_X)^{\ts
m})=0$ {for every} $m\geq 1$ can imply $H^0(X,K_X^{\ts \ell})=0$,
i.e.  $\kappa(X)=-\infty$. By assuming Conjecture \ref{uniruledc},
we obtain Conjecture \ref{geo1}.

Conjecture \ref{mainconjecture1} $\Longrightarrow$ Conjecture
\ref{geo1}. It follows from Theorem \ref{AA} and the fact that
rationally connected manifolds are uniruled.

Conjecture \ref{geo1} $\Longrightarrow$ Conjecture
\ref{mainconjecture1}. The proof follows from Theorem \ref{AA} and a
well-known argument in algebraic geometry (e.g.
\cite[Theorem~1.1]{CDP}, \cite[Corollary~1.7]{GHS03},
\cite[Proposition~2.1]{LP17}), which is also very similar to that of
Theorem \ref{Z0}. Suppose
$$  H^0(X,(T^*_X)^{\ts m})=0, \qtq{for every} m\geq 1.$$
By assuming Conjecture \ref{geo1}, we know $K_X$ is not
pseudo-effective. Hence $X$ is uniruled, thanks to the classical
result of \cite{BDPP13}. Let $\pi:X\dashrightarrow Z$ be the
associated MRC fibration of $X$. After possibly resolving the
singularities of $\pi$ and $Z$, we may assume that $\pi$ is a proper
morphism and $Z$ is smooth. By \cite[Corollary~1.4]{GHS03}, it
follows that the target of the MRC fibration is either a point or a
positive dimensional variety which is not unruled. Suppose $X$ is
not rationally connected, then $\dim Z\geq 1$. Hence $Z$ is not
uniruled, and by \cite{BDPP13} again, $K_Z$ is pseudo-effective. By
assuming Conjecture \ref{geo1}, we obtain
$$H^0(Z, (T^*_Z)^{\ts m_0})\neq 0$$
for some positive integer $m_0$. We obtain $ H^0(X,(T^*_X)^{\ts
m_0})\neq 0$ since $(T^*_Z)^{\ts m_0} \subset (T^*_X)^{\ts m_0}$.
This is a contradiction. \eproof

 If $X$ is rationally connected, $K_X$ is not pseudo-effective.
 Hence there exists a Hermitian metric  $h$ on $K_X^{-1}=\Lambda^{\dim X}T_X$ such
 that $(K_X^{-1}, h)$ is RC-positive. We propose a
 generalization of this fact:

\begin{question}\label{88} Let $X$ be a rationally connected projective manifold.
 Do
there exist  smooth Hermitian metrics $h_p$ on  vector bundles
$\Lambda^p T_X$ ($1\leq p\leq \dim X$) such that $(\Lambda^p T_X,
h_p)$ are all RC-positive? Do there exist  smooth Hermitian metrics
$g_p$ on  vector bundles $ T^{\ts p}_X$ ($p\geq 1$) such that
$(T^{\ts p}_X, g_p)$ are all RC-positive?
\end{question}

\noindent A natural generalization of Conjecture \ref{geo1} is

\begin{question} Let $X$ be a projective manifold.
Suppose $$  H^0(X,(T^*_X)^{\ts m})=0, \qtq{for every} m\geq 1.$$  We
can ask the same question as in Question \ref{88}.
\end{question}

\noindent \textbf{7.2. A partial converse to the Andreotti-Grauert
theorem: the vector bundle version.} We propose a question on vector
bundles converse to Theorem \ref{A}:

\begin{question}\label{7.8} Let $X$ be a projective manifold and $E$ be a vector bundle. Suppose  for every vector bundle $A$,
there exists a positive integer $c_A=c(A,E)$  such that \beq
H^0(X,\mathrm{Sym}^{\ts \ell}E^*\ts A^{\ts k})=0\label{VSYM}\eeq for
$\ell\geq c_A(k+1)$ and $k\geq 0$.  Do there exist smooth Hermitian
metrics $h_p$ on vector bundles $\Lambda^p E$ ($1\leq p\leq
\mathrm{rank}(E)$) such that $(\Lambda^p E, h_p)$ are all
RC-positive? Do there exist smooth Hermitian metrics $g_p$ on vector
bundles $ E^{\ts p}$ (resp. $\mathrm{Sym}^{\ts p}E$) ($p\geq 1$)
such that $(E^{\ts p}, h_p)$ (resp. $(\mathrm{Sym}^{\ts p}E,h_p)$)
are all RC-positive?

\end{question}

\noindent \textbf{7.3. Existence of RC-positive metrics on  vector
bundles} In this subsection, we propose several questions on the
existence of RC-positive metrics.  The celebrated Kodaira embedding
theorem establishes that a line bundle is ample if and only if
carries a smooth metric with positive curvature. The analogous
correspondence for vector bundles is proposed by P. Griffiths:

\begin{conjecture}[\cite{Gri69}] \label{Griffiths1} If $E$ is an ample vector bundle over a compact complex manifold $X$, then $E$ admits a Griffiths positive Hermitian metric.
\end{conjecture}
\noindent When $\dim X=1$, this conjecture is proved in \cite{CF90}.
The following conjecture can be implied by Griffiths' Conjecture
\ref{Griffiths1}:

\begin{conjecture} If $E$ is an ample vector bundle over a projective manifold $X$, then there exists a smooth Hermitian metric $h$ on $E$ such that \bd
\item[(1)] $(E^{\ts k}, h^{\ts k})$ is RC-positive for every $k\geq
1$;
\item[(2)] $(\Lambda^pE,\Lambda^p h)$ is
RC-positive for every $1\leq p\leq \mathrm{rank}(E)$.
\ed\end{conjecture}

\noindent As a converse to Proposition \ref{G}, we also propose the
following
\begin{question}\label{con} Let $E$ be a vector bundle over a compact complex
manifold $X$. If $E$ is weakly RC-positive,  is $E$ necessarily
RC-positive?
\end{question}

\noindent When $\dim X=1$, Question \ref{con} has an affirmative
answer, thanks to \cite{CF90}.

\vskip 1\baselineskip

\section*{Acknowledgements}
I am very grateful to Professor
 Kefeng Liu for his support, encouragement and stimulating
discussions over  years. I would also like to thank Professors
Junyan Cao, J.-P. Demailly, Jixiang Fu, Conan Leung, Si Li, Xiaonan
Ma, Xinan Ma, S. Matsumura, Yum-Tong Siu, Xiaotao Sun, Valentino
Tosatti, Jian Xiao, Weiping Zhang, Fangyang Zheng and Xiangyu Zhou
for helpful suggestions. Finally, I wish to express my sincere
gratitude to Professor S. T. Yau for sharing his deep insights on
the existence of rational curves in the differential geometric
setting and his support.

\end{document}